\documentclass[11pt]{amsart}
\usepackage{geometry}                % See geometry.pdf to learn the layout options. There are lots.
\usepackage{graphicx}
\usepackage{tikz}
\usepackage{pst-node}%,pst-text,pst-3d,pstricks}
\usepackage{paralist}
\usepackage{hyperref}

\usepackage{amssymb}
\usepackage{epstopdf, xyling}
\newtheorem{thm}{Theorem}[section]
\newtheorem{prop}[thm]{Proposition}
\newtheorem{lem}[thm]{Lemma}
\newtheorem{cor}[thm]{Corollary}

\newtheorem{qst}[thm]{Question}
\newtheorem*{thm*}{Theorem}

\theoremstyle{definition}
\newtheorem{dfn}[thm]{Definition}
\newtheorem{eg}[thm]{Example}

\theoremstyle{remark}
\newtheorem{rmk}[thm]{Remark}

\renewcommand{\phi}{\varphi}

\newcommand{\coker}{\mathop{\mathrm{Coker}}\nolimits}

\newcommand{\iso}{\cong}
\renewcommand{\bar}{\protect\overline}

\renewcommand{\i}[1]{\mathfrak{#1}}
\newcommand{\p}{\i{p}}

\newcommand{\rank}{\mathop{\mathrm{rank}}\nolimits}

\newcommand{\depth}{\mathop{\mathrm{depth}}\nolimits}

\newcommand{\h}{\mathop{\mathrm{ht}}\nolimits}

\newcommand{\spec}{\mathop{\mathrm{Spec}}\nolimits}

\newcommand{\pd}{\mathop{\mathrm{pd}}\nolimits}

\renewcommand{\h}{\mathop{\mathrm{ht}}\nolimits}

\newcommand{\lto}{\mathop{\longrightarrow\,}\limits}

\title[Projective Dimension of codimension two algebras presented by quadrics]{The Projective Dimension of codimension two algebras presented by quadrics}
\author{Craig Huneke}
\address{University of Virginia,  Department of Mathematics, 141 Cabell Drive, Kerchof Hall, Charlottesville, VA 22904}
\email{huneke@virginia.edu}
\author{Paolo Mantero}
\address{University of California Riverside,  Department of Mathematics, 900 University Ave., Riverside, CA 92521}
\email{mantero@math.ucr.edu}
\author{Jason McCullough}
\address{Rider University, Department of Mathematics, 2083 Lawrenceville Road,
Lawrenceville, NJ 08648}
\email{jmccullough@rider.edu}
\author{Alexandra Seceleanu}
\address{University of Nebraska-Lincoln,  Department of Mathematics, 203 Avery Hall, Lincoln, NE 68588}
\email{aseceleanu2@math.unl.edu}

\begin{document}
\maketitle

\begin{abstract}
We prove a sharp upper bound  for the projective dimension of ideals of height two generated by quadrics in a polynomial ring with  arbitrary large number of variables. 
\end{abstract}

\section{Introduction}

%Algebras defined as quotients of a polynomial ring by an ideal generated by homogeneous quadratic polynomials have 
%been intensely studied from various perspectives. One reason why they have attracted attention is because they arise as 
%homogeneous coordinate rings of smooth projective varieties embedded by a sufficiently positive ample divisor (\cite{EL}). 
%Another reason is that the class of algebras with defining ideals generated by quadratic polynomials  contains theF
%well-studied class of Koszul algebras (\cite{Froberg}).  For other perspectives on the study of quadratic algebras see \cite{MN}, .... 

This paper is concerned with bounding the projective dimension of ideals generated by homogeneous quadratic polynomials. 
Our main motivation arises from from the following question posed by Stillman:

\begin{qst}[{Stillman \cite[Problem 15.8]{PS}}]\label{Qstillman} Is there a bound, independent of N, on the projective dimension of ideals in $R = K[x_1, . . . , x_N]$ which are generated by $n$ homogeneous
polynomials of given degrees $d_1,\ldots, d_n$?
\end{qst}

Recently, there has been a surge of interest in Stillman's question. We mention here some of the relevant works in order of their appearance: \cite{EnghetaT}, \cite{Engheta1}, \cite{Engheta3}, \cite{McCullough}, \cite{BMNSSS}, \cite{AH}. We also remark that this question has interesting connections to the study of Castelnuovo-Mumford regularity. As shown by Caviglia \cite{Caviglia}, Question~\ref{Qstillman} is equivalent to a similar problem in which  projective dimension is replaced by regularity. For a survey of the developments regarding Stillman's question up to the moment when this paper was written see \cite{MS}.

Although Stillman's question is open in general, it has been answered in the affirmative by Ananyan and Hochster \cite{AH} in the case of ideals generated  by not necessarily homogeneous polynomials of degree at most two. These authors have announced upcoming improvements of their estimates as well as an extension of their result to ideals generated by cubics in \cite{Hoch}. However, the bounds on projective dimension produced by the methods in \cite{AH} are typically very large. Ananyan-Hochster find an upper bound with asymptotical order of growth $2n^{2n}$ for the projective dimension of an ideal generated by $n$ polynomials of degree at most 2 (\cite[Section 6]{AH}). When the number of minimal generators for the ideal is small, the  bounds obtained by Ananyan-Hochster are still large. For example, in the case of ideals generated by three polynomials of degree two, the concrete bound obtained by using the methods in \cite{AH} is 296 \cite[Proposition 3.15]{MS}, while the tight upper bound is 4 by a result of Eisenbud-Huneke (cf. \cite[Theorem 3.1]{MS}). %To see that the bound of 4 is attained consider the ideal $I=(x^2,y^2,ax+by)\subseteq K[a,b,x,y]$.

Our paper stems from a desire to provide a {\em sharp} upper bound for the projective dimension of ideals generated by  homogeneous quadratic polynomials.  This is currently out of reach  in complete generality. Our main result gives a complete answer in codimension two:

\begin{thm*}
 For any ideal $I$ of height two generated by $n$ homogeneous quadratic polynomials in a polynomial ring $R$, $\pd (R/I) \leq 2n-2$. Moreover, this bound is tight.
\end{thm*}

%We proceed by classifying the minimal associated primes of the ideals appearing in the statement of our main theorem and deducing upper bounds for the projective dimension in each of . The most difficult case to analyze turns out to be the case of exactly one  linear associated prime. In that context, we introduce a matrix-theoretic approach which connects our ideals to determinantal ideals of matrices of linear forms. The main  

Our  present work  leads to the 
natural question of whether it is possible to find
% a tight bound can be given for the projective dimension of  ideals generated by  quadratic polynomials, 
%without any additional assumptions on their height. This task turns out to be difficult in general, however forthcoming 
%work in \cite{HMMS3} sheds some light on the case of ideals minimally generated by four quadrics. 
%Alternatively, one may ask for 
bounds on the projective dimension of  ideals generated by  quadratic polynomials of a fixed arbitrary height. In section 6 we pose a more specific question on whether a bound on the projective dimension of ideals generated by  quadratic polynomials can be given in terms of the minimal number of generators and the height of the ideal.

This paper is organized as follows: section 2 covers the background needed for our developments, in section 3 we give an outline of the proof of the main result and we prove this result in the simplest cases, sections 4 and 5  develop the material  needed to fill in the details of the proof of the main result in the most difficult case and section 6 presents some open questions which stem from computations motivated by our work.

\section{Background}

\subsection{Preliminaries}

For the rest of this paper, $R$ will denote a standard graded polynomial ring over a field $K$.  
We further let $I$ denote an ideal of $R$ generated by $n$ homogeneous quadratic polynomials. 
We henceforth use the terminology {\em quadric} for  homogeneous quadratic polynomial and {\em ideal of quadrics} for an ideal generated by quadrics.
 
 For a finitely generated graded $R$-module $M$, there exists a unique, up to isomorphism, minimal graded free resolution
$0 \lto F_p \stackrel{d_p}{\lto}  \cdots \stackrel{d_2}{\lto} F_1  \stackrel{d_1}{\lto} F_0 ,$
that is, an exact sequence of finitely generated graded free $R$-modules $F_i$ with $M \iso \coker(d_1)$. The resolution is called \textit{minimal} if $Im(d_i)\subseteq \mathfrak{m}F_{i-1}$ for $i\geq 1$. Here $\mathfrak{m}$ denotes the homogeneous maximal ideal of $R$. The length $p$ of the minimal free resolution of $M$ is called the \textit{projective dimension} of $M$ and is denote $\pd(M) = p$. 

For the case of cyclic modules, computing $\pd(R/I)$ is equivalent to computing $\pd(I)$ as the two are related by $\pd(R/I)=\pd(I)+1$. For uniformity and consistency with \cite{AH}, \cite{BMNSSS}, \cite{Engheta1}, \cite{Engheta3}, \cite{McCullough}, \cite{MS}  we choose to bound the projective dimension of  $R/I$ throughout the paper. We often make use of the next two well-known remarks which relate the projective dimension of an ideal to the projective dimension of its colon ideals and allow for computations of such colon ideals, respectively.

\begin{rmk}\label{pdlemma}
Let $I$ be an ideal in a polynomial ring $R$ and let $\ell$ be any element of $R$.
Then 
\begin{enumerate}
\item $\pd (R/I) \leq \max \{ \pd (R/(I:(\ell))), \pd (R/(I+(\ell)))\}$
\item $\pd (R/(I:(\ell))) \leq \max \{ \pd (R/I), \pd (R/(I+(\ell))) -1\}$
\item $\pd (R/(I+(\ell))) \leq \max \{ \pd (R/(I:(\ell)))+1, \pd (R/I) \}$
\end{enumerate}
\end{rmk}

\begin{rmk}\label{colon}
Let $I$ be an ideal in a polynomial ring $R$ and let $x,f$ be elements of $R$. Then
$$(I+(xf)):(x)=I:(x) +(f).$$
\end{rmk}

%\begin{proof}
%The non-trivial containment to show is  $(I+(xf)):(x)\subseteq I:(x) +(f)$. We proceed towards this goal by picking $g\in (I+(xf)):(x)$. Then
%$$
%\begin{matrix}
%gx \in I+(xf) & \Longleftrightarrow & gx= h+r xf  & \mbox{ for some } h\in I, r\in R \\
    %                  & \Longleftrightarrow  & h=x(g-rf)  \\
        %              & \Longleftrightarrow  & g-rf \in I:(x)  \\
             %         & \Longleftrightarrow  &g\in  I:(x) +(f).\\
 %\end{matrix}
% $$
%\end{proof}

%KEEP ABOVE PROOF OR DELETE?

\subsection{ 1-generic matrices and heights of ideals of minors}

The class of $1$-generic matrices was introduced in Eisenbud's paper \cite{Eisenbud3}.  We explain in Section 4 how matrices of linear forms come up naturally when analyzing ideals generated by quadrics that are contained in a linear prime.

Let $M$ be a matrix of linear forms over a polynomial ring
$R$ with coefficient field $K$. By a  {\em generalized row} of $M$ we mean a
$K$-linear combination of the rows of $M$ with not all coefficients 0. Similarly a
{\em generalized column} is a nonzero $K$-linear combination of the columns of $M$.
A {\em generalized entry} of $M$ is simply a linear combination with nonzero 
coefficients of the entries of some generalized row.  In the following we write $I_i(M)$ for the ideal  of $i\times i$ minors of $M$.

%\begin{dfn}
%Let $M$ be a $p\times q$ matrix and let $0\leq i\leq \min\{p,q\}$ be an arbitrary integer.%\end{dfn}

\begin{dfn}
We call $M$  {\em $1$-generic} if, after arbitrary $K$-linear row and column operations, $M$ exhibits  no  generalized zero entries. \end{dfn}

 The following result was established in \cite{Eisenbud2} (see also Theorem~2.1 in \cite{Eisenbud3} for a generalization allowing linear sections of small codimension of these determinantal varieties). 

\begin{thm}[Eisenbud, {\cite[Theorem~6.4]{Eisenbud2}}]\label{Eis}
If $M$ is a $1$-generic matrix of linear forms of size $p \times q$ ($p \leq q$) 
with entries in a polynomial ring $R$ over an algebraically closed field $K$, then the
ideal $I_p(M)$ generated by the maximal minors of $M$ is prime of codimension
$q - p + 1$. Its free resolution is given by an Eagon-Northcott
complex and $R/I_p(M)$ is a Cohen-Macaulay domain.
\end{thm}

In addition to using the preceding theorem for 1-generic matrices, we shall also be concerned with ideals of minors of matrices which are far from being 1-generic; specifically we shall be interested in determining the height of ideals of minors of matrices whose rank is not maximal. In this 
situation, bounds on the height of ideals of minors for (not necessarily linear) matrices have been given by Eisenbud-Huneke-Ulrich \cite{EHU}
generalizing results of Eagon-Northcott, Bruns and Faltings. 

\begin{thm}[Eisenbud-Huneke-Ulrich, {\cite[Theorem~A]{EHU}} ]\label{EHU}
 Let $R$ be a regular local ring, and let $M$ be a matrix of size 
$p\times q$ with entries in $R$. Set $r=\rank(M)$ and 
consider an integer $i$ such that $1\leq i \leq r$ and $I_i (M)\neq R$. Then
$$\h I_i(M) \leq  (r - i + 1)( \max\{p, q\}-i + 1) + i - 1.$$
\end{thm}

\begin{rmk}
The conclusion of the theorem above holds in the graded case as well.  If $R$ is a polynomial ring and $M$ a matrix of size 
$p\times q$ and rank $r$ whose entries are homogeneous forms, then
$$\h I_i(M) \leq  (r - i + 1)( \max\{p, q\}-i + 1) + i - 1.$$
This is easily deduced from theorem \ref{EHU} by localizing at the homogeneous maximal ideal.
\end{rmk}

\subsection{Residual intersections}

 The second result that we use extensively appeared in the context of residual intersections, a notion introduced by Artin and Nagata \cite{AN}. An ideal $J$ is called an $s$-{\it residual intersection} of an ideal $I$ if there exists an $s$-generated ideal $A\subseteq I$ such that $J=A:I$ and $h(J)\geq s\geq\h(I)$. If, in addition, $\h (I+J)> \h (J)$, then $J$ is said to be an $s$-{\it geometric residual intersection} of $I$. The notion of residual intersection generalizes that of liaison, for which $s=\h(I)$. 
 
One of the most common settings in which residual intersections arise naturally is the following: take $A$ to be an ideal generated by at most $s$ elements, assume $A$ is not unmixed, and let $I$ be the intersection of its primary components of height at most $s$ (in any of its irredundant primary decompositions). If $I\neq A$, then $J=A:I$ is a geometric residual intersection of $I$, in particular $J$ is actually the intersection of all the primary components of $A$ of height at least $s+1$ and  $A=I\cap J$.
 
%The main motivation for Artin and Nagata to introduce this notion is the following question: given an ideal $A=I\cap J$, what properties of $I$ and $A$ imply that $J$ is `nice'?
% In \cite{Huneke}, the first author proved the following result, which we rely on in our proofs: if $I$ is strongly Cohen-Macaulay and satisfies the property $G_s$, then $J=A:I$ is a geometric residual intersection (yielding $A=I\cap J$) and $J$ is Cohen-Macaulay. We begin by defining the strongly Cohen-Macaulay and $G_s$ properties.

\begin{dfn}
Let $I$ be an ideal in  a polynomial ring $R$ and fix $\underline{f}=f_1,\ldots,f_n$ to be any minimal set of generators of $I$. We define $I$ to be strongly Cohen-Macaulay if the  homology modules $H_i(\underline{f},R)$ of the Koszul complex associated to $\underline{f}$ are either zero or Cohen-Macaulay modules for all $i$.
\end{dfn}

Although this definition seems to depend on the chosen generating set $\underline f$ of $I$, one can check that this is actually a property of the ideal $I$ (see for instance \cite[Remark~1.2]{Huneke}).
%As noted in \cite{Huneke}, the most obvious example of strongly Cohen-Macaulay modules are the ones where $I$ is locally a complete intersection.

\begin{dfn}[Artin-Nagata\cite{AN}]
We say $I$ satisfies condition $G_s$ if $\mu(I_{\mathfrak{p}})\leq \dim R_{\mathfrak{p}}$ for every $\mathfrak{p}\in Spec(R)$,  with $I \subset \mathfrak{p}$ and $\dim R_{\mathfrak{p}}\leq s-1$.
\end{dfn}

In due course we shall exploit the following result which is a translation of  \cite[Theorem~3.1]{Huneke} into purely algebraic language.

\begin{thm}$($Huneke, \cite[Theorem~3.1]{Huneke}$)$\label{Craig}
Let  $R$ be a Cohen-Macaulay local ring. Let $I$ be an ideal of $R$ such that $I$ is strongly Cohen-Macaulay and  satisfies $G_s$. Let $J=A:I$ be an $s$-geometric residual intersection of $I$. Then:
\begin{enumerate}
\item $R/J$ is Cohen-Macaulay and $\h J= s$, 
\item $A=I\cap J$,  
\item $I+J$ is strongly Cohen-Macaulay and 
\item $\depth R/A\geq \dim R-s$
\end{enumerate}
\end{thm}

The interested reader may find generalizations of this result by Herzog-Vasconcelos-Villareal and Huneke-Ulrich in \cite{HVV} and \cite{HU}, respectively. However, for the purpose of this paper, we only need the following special case of Theorem \ref{Craig}, where $I$ is taken to be a homogeneous complete intersection ideal. 
\begin{cor}\label{residual}
Let  $R$ be a polynomial ring. Let $C$ be a complete intersection homogeneous ideal in $R$. Let $A=(a_1,\ldots a_s)$ denote an ideal contained in $C$. 
Set  $J =A:C$ and assume $\h J \geq s$ and $\h(C+J)\geq s+1$.  Then $\h (J)=s$ and $\pd (R/A)\leq s$.
\end{cor}

\begin{proof}
The problem can be reduced to the local case by localizing with respect to the homogeneous maximal ideal $\mathfrak{m}$ of $R$. Clearly $I_{\mathfrak{m}}$ is strongly Cohen-Macaulay and $G_s$, hence we can apply Theorem \ref{Craig} to obtain $\h J_{\mathfrak{m}}=s$ and $\depth R_{\mathfrak{m}}/A_{\mathfrak{m}}\geq \dim R_{\mathfrak{m}}-s$. This latter inequality combined with  the Auslander-Buchsbaum formula implies $\pd(R_{\mathfrak{m}}/A_{\mathfrak{m}})\leq s$. The conclusion now follows from the equalities $\h(J)=\h(J_{\mathfrak{m}})$ and $\pd_R (R/A) =\pd_{R_{\mathfrak{m}}} (R_{\mathfrak{m}}/A_{\mathfrak{m}})$.
\end{proof}

\subsection{Hilbert-Samuel multiplicity}

In our quest to analyze the projective dimension of ideals generated by quadrics, we shall resort to a case analysis based on the minimal primes associated to our ideals. To facilitate the understanding of what minimal primes may occur, we  use the notion of Hilbert-Samuel multiplicity and a result known as the associativity formula. In the following,  we denote by $e(M)$ the Hilbert-Samuel multiplicity of a module $M$ and by $\lambda(M)$ the length of an Artinian module $M$.

%For a graded module $M$, the {\it Hilbert series} $H_M(t)=\sum_{i\geq 0}dim_K M_i t^i$ can be written as a rational function of the form $H_M(t)=\frac{h(t)}{(1-t)^d}$, where $d$ is the dimension of $M$ and $h$ is a polynomial of degree at most $N$. We define the {\it multiplicity} of a graded $R$-module $M$ to be the value $e(M)=h(1)$. For an artinian module $M$, the multiplicity is equal to the {\it length} of the module defined as $\lambda(M)=\sum_{i\geq 0} dim_K M_i$. By convention, for a homogeneous ideal $I$, we refer to $e(R/I)$ as the multiplicity of $I$.

\begin{thm}[Associativity Formula {\cite[Theorem 14.7]{Matsumura}}]\label{assoc}  Let $I$ be an ideal of $R$.  Then
\[e(R/I) = \sum_{\substack{\p \in \spec(R)\\\h(\p) = \h(I)}} e(R/\p) \lambda(R_\p/J_\p).\]
\end{thm}

Let $I^{un}$ denote the unmixed part of $I$, defined as the intersection of the primary components of $I$ with height equal to $\h(I)$. It is easy to see from the formula above that $e(R/I)=e(R/I^{un})$.

 We recall a classical lower bound for the multiplicity of non-degenerate prime ideals. A homogeneous ideal is called degenerate if it has a linear form as a  minimal generator.

\begin{prop}( \cite[Proposition~0]{EH2}, \cite[Corollary~18.12]{EH1})\label{non-deg}
Let $\mathfrak{p}$ be a homogeneous prime ideal in a polynomial ring defined over an algebraically closed field. If $\mathfrak{p}$ is non degenerate, then $e(R/\mathfrak{p})\geq \h (\mathfrak{p})+1$.
\end{prop}

In view of this bound,  an  ideal $\mathfrak{p}$ is said to have {\em minimal multiplicity} if $e(R/\mathfrak{p})=\h (\mathfrak{p})+1$. A variety defined by a 
prime ideal of minimal multiplicity is called a variety of minimal multiplicity. These varieties are well known and have been classified. Proofs of this classification in characteristic zero can be found in \cite{Bertini}, \cite{Harris}, \cite{Xambo} and a characteristic-independent proof can be found in \cite{EH2}. The case of interest to us, that of varieties of  codimension 2 and minimal multiplicity, has been first considered in \cite{X}, \cite{Swinnerton-Dyer}.  

\begin{thm}[Del Pezzo-Bertini, {\cite[Theorem~1]{EH2}}]\label{minmultprimes}
A non-degenerate irreducible variety of codimension at least two and minimal multiplicity  defined over an algebraically closed field is of one of the following types
\begin{inparaenum}
\item[(1)] a rational normal scroll,
\item[(2)] the second Veronese embedding of $\mathbb{P}^2$ in $\mathbb{P}^5$ or  
\item[(3)] a cone over one of the previous two varieties.
\end{inparaenum}
\end{thm}

\section{Main results}

\subsection{The method of proof}

The main result of this paper provides a tight upper bound on the projective dimension of ideals $I$ of height two generated by quadrics. The first step in our proof is to classify  the minimal primes  of such ideals. In order to accomplish this task, we  bound the multiplicity of height two ideals generated by quadrics  then we use the associativity formula to find the possible multiplicities of the individual minimal primes. As a consequence of this analysis, we conclude that a height two ideal of quadrics $I$ is contained in at least one prime $\mathfrak{p}$ of  one of the following types: 
\begin{inparaenum}
\item[(1)] a prime of multiplicity one and height two, i.e.  $\mathfrak{p}=(x,y)$, with $x,y$ independent linear forms,
\item[(2)] a prime of multiplicity two and height two, i.e.  $\mathfrak{p}=(x,q)$, with $x$ a linear form and $q$ an irreducible quadic  or 
\item[(3)] a prime of multiplicity three and height two, i.e the defining ideal of one of the varieties of minimal multiplicity classified in Theorem \ref{minmultprimes} (we shall also see that degenerate ideals of multiplicity three and height two cannot occur).
\end{inparaenum}

We prove bounds on the projective dimension of ideals of quadrics contained in primes of each of the types described above separately. The first case is the most intricate and is dealt with in Section 5, while the other two cases are covered in this section. We combine these bounds in Theorem \ref{main}, which conveys the overall estimate. 

\subsection{Minimal associated primes  and consequences on projective dimension}

We begin by characterizing the minimal primes associated to height two ideals of quadrics and classifying them according to their multiplicity.

\begin{lem}\label{Pmult} Let $I$ be an ideal of $R$ minimally generated by $n\geq 3$ quadrics and with $\h(I) = 2$, then \mbox{$e(R/I) \le 3$}.
\end{lem}

\begin{proof}
Let $f, g$ be a regular sequence of quadratic forms contained in $I$, thus $e(R/(f, g)) = 4$.  Passing to the artinian reduction of $R/I$ by going modulo a sequence of linear forms $\underline{\ell}$, yields $e(R/I)\leq\lambda(R/(\underline{\ell},I)) \leq \lambda(R/(\underline{\ell},f,g))=e(R/(f,g))= 4$. We have the series of containments $(f,g) \subset I \subset I^{un}$.  Note that $(f,g)$ and $I^{un}$ are unmixed ideals of height two.  If $e(R/I)=e(R/I^{un})=e(R/(f,g))=4$, then $(f,g) = I^{un}$ by \cite[Lemma 8]{Engheta1}.  But this would force $(f,g) = I$, contradicting that $I$ has at least  three minimal generators.  Thus we must have $ e(R/I)\leq 3$.
\end{proof}

Combining this lemma with the associativity formula one obtains:

\begin{cor}\label{primes}
Let $I$ be an ideal of $R$ minimally generated by $n\geq 3$ quadrics and with $\h(I) = 2$. Let $\mathfrak{p}$ be an minimal prime of $I$ with $\h \mathfrak{p}=2$. Then $e(R/\mathfrak{p})\leq 3$.
\end{cor}

To make use of this result, one needs to understand prime ideals of height two and multiplicity 1, 2 and 3. According to Theorem \ref{minmult}, any prime of height two and multiplicity 1 or 2 is degenerate (contains a linear form). Going modulo this linear form, the image is a principal ideal. In particular, a prime 
of height two and multiplicity one is minimally generated by two independent linear forms, and a prime 
of height two and multiplicity two is minimally generated by a linear form and a quadric.  This gives the first two types of primes mentioned in the paragraph at the beginning of this section. Similarly, it is easy to see that the degenerate primes of height two and multiplicity three are minimally generated by a linear form and a cubic, while the non-degenerate ones are listed in Theorem \ref{minmultprimes}.

We begin our case analysis by considering  ideals generated by quadrics contained in a prime ideal of multiplicity three and height two.

\begin{prop}\label{minmult}
Let $R$ be a polynomial ring  over an algebraically closed field. Any height two ideal  $I$ generated by quadrics contained in a prime ideal $\mathfrak{p}$ with $\h (\mathfrak{p})=2$ and $e(R/\mathfrak{p})=3$ has $\pd (R/I )= 2$ . In particular, for any such ideal $I$, $R/I$ is Cohen-Macaulay.
\end{prop}

\begin{proof} 
We start by considering the degenerate case $\mathfrak{p}=(x,c)$, with $x$ a linear form, $c$ an irreducible cubic polynomial and $\h(x,c)=2$. Since $I\subseteq \mathfrak{p}$ is generated by quadrics,  all minimal generators of $I$ must be multiples of $x$. However, in that case the height of $I$ is one, a contradiction.

Next assume $\mathfrak{p}$ is  the non-degenerate defining ideal of one of the varieties listed in Theorem \ref{minmultprimes}. The ideal of the second Veronese of $\mathbb{P}^2$ has height three and must be excluded. Note that each of the defining ideals of rational normal scrolls of height two is minimally generated by 3 quadrics. Since the minimal generators of $I$ are linear combinations of these quadrics, we obtain, depending on the number of minimal generators of $I$, that either $I$ is a complete intersection or $I=\mathfrak{p}$. In the first case $\pd (R/I)=2$ and in the second case $\pd (R/I)=2$ as well, as $R/\mathfrak{p}$ is Cohen-Macaulay by the Hilbert-Burch theorem applied to each of the primes defining rational normal scrolls.
\end{proof}

We continue with the case of ideals generated by quadrics contained in a prime ideal of multiplicity  two and height two.

\begin{prop}\label{xq}
Let $I$ be an ideal  generated by $n$ quadrics  which is contained in a prime ideal of height two and multiplicity two. Then $\pd (R/I) \leq n$.
\end{prop}

\begin{proof}
Let $\mathfrak{p}$ be the minimal prime mentioned in the hypothesis. It is necessarily of the form $\mathfrak{p}=(x,q)$, with $x$ a linear form and $q$ a quadric.
We may assume that:
$$I=(q+x\ell_1, x\ell_2,\ldots,x\ell_n)$$
where $\ell_i$ are linear forms such that $\{\ell_2, \ldots, \ell_n\}$ form a linearly independent set. By first using Remark \ref{colon} and subsequently using that $x$ is regular on $R/(q)$ and hence also on $R/(q+x\ell_1)$, we compute
$$I:(x)=(q+x\ell_1):(x)+( \ell_2,\ldots,\ell_n)= (q+x\ell_1)+( \ell_2,\ldots,\ell_n).$$
If $q+x\ell_1 \in (\ell_2,\ldots,\ell_n)$, then $I:(x)=(\ell_2,\ldots,\ell_n)$ and $\pd (R/I:(x)) = n-1$. Otherwise, $q+x\ell_1$ is a regular element on $R/(\ell_2,\ldots,\ell_n)$ and thus 
 $\pd (R/(q+x\ell_1, \ell_2,\ldots,\ell_n)) =n$. Hence in any case we obtain $\pd (R/(I:(x)))\leq n$.
Now   since $I+(x)=(q,x)$, we obtain $\pd (R/(I+(x)))=2$ and  Remark \ref{pdlemma} yields $\pd (R/I)\leq n$.
\end{proof}

The last and most involved case is that of ideals generated by quadrics contained in a  prime ideal of multiplicity one. We defer the detailed analysis of this case to Section 5. %mentioning only that as a main result we obtain in Proposition \ref{linear} a bound of $2n-2$ on the projective dimension of cyclic modules defined by ideals generated by $n$ quadrics.

\subsection{Proof of the main theorem}

The following theorem, which is our main result, combines the bounds obtained throughout this paper.

\begin{thm}\label{main}
 For any ideal $I$ of height two generated by $n$ quadrics in a polynomial ring $R$, one has $\pd (R/I) \leq 2n-2$. 
\end{thm}

\begin{proof}
  For the purpose of computing projective dimension we may assume the ground field is algebraically closed by tensoring $R$ with an appropriate extension of the original field. The case when $n=2$ yields $R/I$ is a complete intersection, thus in this case $\pd (R/I) =2$  agrees with the conclusion of the theorem. Henceforth we assume $n\geq 3$.
 
   By Corollary \ref{primes}, any ideal of height two generated by quadrics falls under the hypotheses of at least one of Proposition \ref{linear}, Proposition \ref{xq} or Proposition \ref{minmult}. These three results insure that $\pd (R/I) \leq 2n-2$,  $\pd (R/I) \leq n+1$ or $\pd (R/I )=2$ respectively. 
 Taking the maximum of these bounds yields $\pd (R/I)\leq 2n-2$ overall.
\end{proof}

We complement the statement of our principal result by the important observation that the bound in Theorem~\ref{main} is tight as we demonstrate in Example~\ref{tight}.  The example below is a particular case of the main result in \cite{McCullough}, specifically the ideal discussed below appears under the name $I_{2,n-2,2}$ in \cite{McCullough}. Note that the multiplicity of the family of examples  given below is one as soon as the minimal number of generators is at least 
four and it is two for the example minimally generated by three quadrics.

\begin{eg}\label{tight}
Let $n\geq 2$ and $R=K[x,y,a_{1,1},\ldots,a_{2,n-2}]$. Consider the $n$-generated ideal of quadrics  
$I=(x^2,y^2,a_{1,1}x+a_{2,1}y,\ldots , a_{1,n-2}x+a_{2,n-2}y).$
Then $\pd (R/I)= 2n-2$.
\end{eg}

\begin{proof}
 By the graded Auslander-Buchsbaum formula, showing that $\pd (R/I)= 2n-2$ is equivalent to showing $\depth R/I=0$ or, equivalently, that the maximal ideal of $R$ is associated to $R/I$. It is not hard to see that  $xy$ is a non-zero socle element in $R/I$, which finishes the proof. 
 \end{proof}
 
\subsection{A more detailed analysis}
Even though our bound is tight, it turns out that it can be refined in an important number of cases. The purpose of the additional analysis we perform in the rest of this section is to point out in which contexts such refinements are possible. Following Engheta \cite{EnghetaT}, we introduce the following notation to keep track of the possibilities for the associated primes of minimal height of an ideal $J$.  
\begin{dfn}\label{type}
We say $I$ is of type $\langle e = e_1, e_2, \ldots, e_m ; \lambda = \lambda_1, \lambda_2,\ldots, \lambda_m\rangle$ if $I$ has $m$ associated primes $\p_1,\ldots, \p_m$ of minimal height with $e(R/\p_i) = e_i$ and with $\lambda(R_{\p_i}/I_{\p_i}) = \lambda_i$, for $i = 1,\ldots, m$.  
\end{dfn}

\begin{rmk}
If $I$ is of type $\langle e = e_1, e_2, \ldots, e_m ; \lambda = \lambda_1, \lambda_2,\ldots, \lambda_m\rangle$, then 
$e(R/I)=\sum_{i=1}^m e_i\lambda_i$. 
\end{rmk}
We classify the possible types of height two ideals generated by quadrics as follows:

\begin{prop}\label{list}
Let $I$ be an ideal of $R$ minimally generated by $n\geq 3$ quadrics and with $\h(I) = 2$, then the minimal associated primes of $I$ fall into one of the following categories
\begin{inparaenum}
\item[(1)] $\langle1;1 \rangle$
\item[(2)] $\langle2;1 \rangle$
\item[(3)] $\langle1;2 \rangle$ 
\item[(4)] $\langle1,1;1,1\rangle $
\item[(5)] $\langle 3;1 \rangle$ 
\item[(6)] $\langle 1;3\rangle $ 
\item[(7)] $\langle 1,2;1,1\rangle $ 
\item[(8)] $\langle 1,1;1,2 \rangle$
\item[(9)] $\langle 1,1,1;1,1,1\rangle $ 
\end{inparaenum}.
\end{prop}

\begin{proof}
The  result of combining proposition \ref{minmultprimes} and the associativity formula \ref{assoc} is the inequality $$e(R/I) = \sum_{\substack{\p \in \spec(R)\\\h(\p) = 2}} e(R/\p) \lambda(R_\p/I_\p)\leq 3.$$ The possible types listed above are all the possibilities of partitioning the integers 1,2 and 3  according to the associativity formula.
\end{proof}

A finer case analysis than the one performed in the proof of Theorem \ref{main} is undertaken in the  table at the end of this section, which  summarizes the bounds we are able to obtain on the projective dimension of  ideals of quadrics of each of the types described above. We base our refined analysis on the observation that better bounds can be established for ideals contained in multiple structures with linear support (Proposition \ref{1;n}) or in at least two distinct linear primes (Proposition \ref{2lin}). For each bound, we reference the relevant results in the paper for the reader's convenience. Note that by taking the maximum among the bounds listed in the table, we recover our general result, Theorem \ref{main}.

The case breakdown in the table  identifies that, for ideals of height two minimally generated by at least four quadrics, there is only one type that can achieve the maximum projective dimension, namely  ideals of multiplicity one. For ideals minimally generated by three quadrics, the situation is quite different: the sharp bound of 4 is never attained by an ideal of multiplicity one (see \cite{HMMS2} for a proof). The bound of $n+2$ on the projective dimension of ideals of height two generated by  quadrics and having multiplicity strictly larger than one may be of separate interest. 

\bigskip

\noindent
\begin{tabular}{|c|c|l|}
\hline
$\langle \underline{e};\underline{\lambda}\rangle$&  Bound on $\pd(R/I)$ & Reference\\
\hline
$\langle1;1\rangle$&$2n-2$& contained in a linear prime (see Prop.~\ref{linear})\\
\hline
$\langle2;1\rangle$& $n$  & contained in a multiplicity two prime (see Prop.~\ref{xq})\\
$\langle1;2\rangle$ & $n+2$  & multiple structure on a linear prime (see Prop.~\ref{1;n})\\
$\langle1,1;1,1\rangle$& $n+1$ &  contained in 2 linear primes (see Prop.~\ref{2lin})\\
\hline
$\langle3;1\rangle$ & 2 &  prime of minimal multiplicity (see Prop.~\ref{minmult})\\
$\langle1;3\rangle$ & $n+2$ & multiple structure on a linear prime (see Prop.~\ref{1;n})\\
$\langle1,2;1,1\rangle$ & $n$ & contained in a multiplicity two prime (see Prop.~\ref{xq}) \\
$\langle1,1;1,2\rangle$ & $n+1$  & contained in 2 linear primes (see Prop.~\ref{2lin})\\
$\langle1,1,1;1,1,1\rangle$ & $n+1$ & contained in 2 linear primes (see Prop.~\ref{2lin})\\
\hline

\end{tabular}

\bigskip

\section{Matrices of linear forms  and ideals of quadrics contained in linear primes}

It is easy to see (\cite[Exercise 7.6]{Hartshorne}) that any prime ideal of multiplicity one of a polynomial ring over an algebraically closed field is generated by linear forms. We shall call such a prime ideal generated by
linear forms a {\em linear prime}. 

We give a matrix-theoretic approach to ideals generated by quadrics  contained in a linear prime by viewing them as subideals of the
determinantal ideal associated to certain matrices of linear forms.  While our ideals are not the full determinantal ideal, we often use residual intersection techniques to recover  determinantal ideals as residuals of our ideals or of subideals thereof (see Proposition \ref{nozero} for a typical example). However, this technique requires 1-genericity assumptions  on the matrix of linear forms, therefore we begin by studying linear algebraic properties of matrices associated to ideals of quadrics.

For the rest of the section, we  consider the case of an ideal $I$ generated by $n$ quadrics and such that $I\subset (x,y)$, where $x$ and $y$ are linearly independent linear forms.
\begin{dfn}\label{defA}
Let $$A=\left( \begin{matrix} a_{11} & \ldots & a_{1n} \\  a_{21} & \ldots & a_{2n}\\\end{matrix}\right ). $$
 %$A$ be a $s\times n$ matrix of linear forms in a polynomial ring $R$. 
Let $I$ be an ideal generated by the $n$ quadrics  which are entries of the $1\times n$ matrix 
 %$\left( \begin{matrix} x_1 & \ldots & x_s \end{matrix}\right )A$, where $x_1, \ldots, x_s$ are linear forms.
 $\left( \begin{matrix} -y & x \end{matrix}\right )A$, where $x, y$ are linear forms.
In this context we say that $I$ is  {\em represented by coefficients} by the matrix $A$.
\end{dfn}

\begin{dfn}\label{defM}
Let $$M=\left( \begin{matrix} x & a_{11} & \ldots & a_{1n} \\ y & a_{21} & \ldots & a_{2n}\\\end{matrix}\right )  $$
 be the matrix obtained by prepending the column vector $\left( \begin{matrix} x & y \end{matrix}\right )^T$
 %$\left( \begin{matrix} x_1 & \ldots & x_s \end{matrix}\right )^T$ 
 to the matrix $A$ described in Definition \ref{defA}.
 Let $I$ be the ideal generated by the subset of the $2$ by $2$ minors of this matrix  which involve the first column of $M$. 
In this context we say that $I$ is {\em represented by minors} by the matrix $M$.
\end{dfn}

\begin{rmk}
The same ideal $I$ is obtained from the two matrix representations described in \ref{defA} and \ref{defM} respectively. 
\end{rmk}

The example below shows that the same ideal $I$ can be represented by several distinct matrices both by coefficients and also  by minors .

\begin{eg}
Both  matrices $A$ and $A'$ displayed below represent by coefficients the ideal $$I=(x^2,xy,ax+by,cx+dy)\subset K[x,y,a,b,c,d].$$ Furthermore, both matrices $M$ and $M'$  represent the same ideal $I$ by minors:
$$A=\left( \begin{matrix} 0 & 0 & -b & -d  \\ x & y & a &c  \\\end{matrix}\right ), \quad  A'=\left( \begin{matrix} 0 & -x & -b & -d  \\ x & 0 & a &c  \\\end{matrix}\right ),$$
$$M=\left( \begin{matrix} x &0 & 0 & -b & -d \\ y & x & y & a &c \\\end{matrix}\right ), \quad  M'=\left( \begin{matrix} x &  0 & -x & -b & -d \\ y & x & 0 & a &c   \\\end{matrix}\right ).$$
\end{eg}

Note that $M'$ in the example above is obtained from $M$ by subtracting the first column from the third column. We use elementary operations of this type  frequently this section. We define a {\em sequence of invertible elementary  operations} applied to a matrix to mean multiplication on either side by invertible matrices with scalar entries. Note that this is closely related to the notions of generalized rows and columns  defined in the  introduction. We say that two matrices are {\em equivalent} if one is obtained from the other by performing a sequence of elementary operations. 

\begin{rmk}\label{rowoper}
The ideals of quadrics represented by coefficients by a pair of equivalent matrices are equal. To have that two ideals of quadrics represented by minors by a pair of equivalent matrices are equal, we must further require that any elementary column operations involved preserves the first column. We refer to matrix operations with this property by the name of {\em ideal-preserving operations}.    \end{rmk}

A fundamental idea of our analysis is to systematically exploit the pattern of generalized zero entries of the matrix representing an ideal of quadrics by minors. We now explain a procedure to reduce  matrices of linear forms using ideal-preserving operations to a short list of canonical forms, which shall be analyzed in detail in section 5.

\begin{lem}\label{putzero}
Let $M$ be a $2\times (n+1)$ matrix of linear forms representing a height two ideal by minors and having at least one generalized zero. Then there is a sequence of ideal-preserving operations producing an equivalent matrix $M'$ of the form $$M'=\begin{pmatrix}x' & 0 & a_{12} & \ldots & a_{1n}\\ y' & \ell & a_{22} & \ldots & a_{2n}\\\end{pmatrix},$$
where all the entries of $M'$ are either linear forms or zero and $\h(x',y')=2$. 
\end{lem}

\begin{proof}
Recall that a generalized zero of $M$ is a  linear combination of the entries of a generalized column of $M$. Since $M$ represents a height two ideal, the entries of the first column of $M$ are linearly independent, thus this generalized column cannot be equal to the first column of $M$. Therefore at least one  of the last $n$ columns of $M$ is involved in the linear combination giving the previously mentioned generalized column. The matrix operation replacing this column of $M$ by the generalized column is an ideal-preserving invertible column operation. 
%Therefore the linear combination yielding the  generalized column must involve at least one of the last $n$ columns of $M$. After performing a permutation of the last $n$ columns (an ideal-preserving operation), we may assume that the second column appears with a non-zero coefficient in the linear combination yielding the generalized column containing the generalized zero. The matrix  operation replacing the second column by the  generalized column  is then an invertible ideal-preserving elementary operation. 
Next an appropriate row operation (which is in turn ideal-preserving) can be performed to yield a zero entry in the upper entry of this column, followed by a permutation of the columns which places the zero entry in the position indicated in $M'$.  \end{proof}

\begin{prop}\label{stdzero}
Let $M$ be a $2\times(n+1)$ matrix of linear forms with $n\geq 2$ representing a height two ideal of quadrics by minors. Then $M$ is equivalent via a sequence of ideal-preserving elementary operations to a matrix $M'$ of one of the following types, where the entries of all matrices below are either zero or linear forms and $\h(x,y)=2$:
\begin{enumerate}
\item $M'$ is 1-generic;
%\item $M'$ has exactly one generalized zero which occurs in the last ;
\item $M'=\begin{pmatrix}x & 0 & a_{12} & \ldots & a_{1n}\\ y & a_{21}& a_{22} & \ldots & a_{2n}\\\end{pmatrix}$, where  $D=\begin{pmatrix}x &   a_{12} & \ldots & a_{1n}\\ y & a_{22} & \ldots & a_{2n}\\\end{pmatrix}$ is 1-generic;
\item $M'=\begin{pmatrix}x & 0 & 0 & a_{13} & \ldots & a_{1n}\\ y &a_{21} & a_{22} & a_{23} & \ldots & a_{2n}\\\end{pmatrix}$, with no additional restrictions;
\item $M'=\begin{pmatrix}x & 0 & a_{12} & a_{13} & \ldots & a_{1n}\\ y & a_{21} & 0 & a_{23} & \ldots & a_{2n}\\\end{pmatrix}$, where  $D=\begin{pmatrix}x &   a_{13} & \ldots & a_{1n}\\ y & a_{23} & \ldots & a_{2n}\\\end{pmatrix}$ is 1-generic;
\item $M'=\begin{pmatrix}x & 0 &a_{12} & a_{13} & a_{14} & \ldots & a_{1n}\\ y & a_{21} & 0  & \lambda a_{13} & a_{24} & \ldots & a_{2n}\\\end{pmatrix}$, where $\lambda$ is a scalar and there are  no additional restrictions.
\end{enumerate}
\end{prop}

\begin{proof} We argue based on the number of generalized zeros of $M$.
If $M$ has no generalized zeros, then we set $M'=M$ which is a matrix of the first type listed. If $M$ has at least one generalized zero, then we use Lemma \ref{putzero} to produce by ideal-preserving operations a matrix $M'=\begin{pmatrix}x' & 0 & a_{12} & \ldots & a_{1n}\\ y' & a_{21} & a_{22} & \ldots & a_{2n}\\\end{pmatrix}$. We set $D=\begin{pmatrix}x' &   a_{12} & \ldots & a_{1n}\\ y' & a_{22} & \ldots & a_{2n}\\\end{pmatrix}$ and analyze the pattern of generalized zeros of the submatrix $D$. If $D$ has no generalized zeros then $M'$ is of the second type listed in the statement.
Next we consider the case when $D$ has at least one generalized zero. Using Lemma \ref{putzero}, $D$ is equivalent to a matrix of the type $D'=\begin{pmatrix}x'' & 0 & a'_{13} & \ldots & a'_{1n}\\ y'' & a_{22}' & a'_{23} & \ldots & a'_{2n}\\ \end{pmatrix}$ via  ideal-preserving elementary row operations. Applying the same operations to the matrix $M'$ yields a result of the form
$$M''=\begin{pmatrix}x'' & \alpha a_{21} &0 & a'_{13} & \ldots & a'_{1n}\\ y'' & \beta a_{21} & a_{22}' & a'_{23} & \ldots & a'_{2n} \\ \end{pmatrix}.$$

Now either $\alpha=0$, which gives a matrix of the third type on our list, or $\alpha\neq 0$ and then subtracting a multiple of the first row  from the second yields an equivalent matrix
$$M''=\begin{pmatrix}x'' & \alpha a_{21} &0 & a'_{13} & \ldots & a'_{1n}\\ y''' & 0 & a_{22}' & a''_{23} & \ldots & a''_{2n} \\ \end{pmatrix}.$$
At this stage, either the submatrix $\begin{pmatrix} x''  & a'_{13} & \ldots & a'_{1n}\\ y'' & a''_{23} & \ldots & a''_{2n} \\ \end{pmatrix}$ has no generalized zeros, which gives the fourth form in the statement of this proposition,  or we may use Lemma \ref{putzero} to find ideal-preserving operations which replace $a'_{13}$ by zero. Acting by these matrix operations on $M''$ yields
$$M'''=\begin{pmatrix}x''' & \alpha' a_{21} & \gamma' a_{22} & 0 & a'''_{14}&  \ldots & a'''_{1n}\\ y''' & \beta' a_{21} & \delta' a_{22} & a'''_{23} & a'''_{14} & \ldots & a'''_{2n} \\ \end{pmatrix}.$$
If $\alpha'$ is zero, then the matrix is of the third type on our list up to a permutation of columns. Otherwise,  by subtracting a $\beta'/\alpha'$ multiple of the first row from the last and permuting the third and fourth columns, we obtain a matrix of the last type on our list. %(with $\lambda=\delta'-\beta'\gamma'/\alpha'$).
 \end{proof}

Finally, we give a result concerning colon ideals. We strengthen the containments presented below to equalities in the case when $A$ is a 1-generic matrix in Proposition \ref{nozero}.

\begin{lem}\label{Cramer}
Let  $I$ be an ideal of quadrics $I$ contained in a height two linear prime $(x,y)$ and let $A$ and $M$ be matrices representing $I$ by coefficients and minors respectively. Then
\begin{enumerate}
%\item $I_s(A)\subseteq I:(x_i)$, for any $1\leq i \leq n$
\item $I_2(A)\subseteq I:(x,y)$
\item $I_2(M) \subseteq I:(x,y)$
\end{enumerate}
\end{lem}

\begin{proof}
This first assertion is a generalization of the familiar Cramer's rule. The second assertion puts together the first containment and the trivial containment $I \subseteq I:(x,y)$, using that $I_2(M)=I_2(A)+I$.
\end{proof}

\section{Bounding the projective dimension of ideals of quadrics contained in a linear prime of height two}

%\subsection{The case of a 1-generic matrix}

We henceforth start an analysis of ideals $I$ generated by $n\geq 2$ quadrics which are contained in a linear prime ideal of height two by considering ideals represented  by each of the matrices $M$ listed in Proposition \ref{stdzero}. Recall that we work with ideals of quadrics in a polynomial ring $R$ over a field $K$. For the rest of this section we assume that $K$ is  algebraically closed   so that we can apply Theorem \ref{Eis} on 1-generic matrices.  This  assumption is easily met by tensoring with an extension of the original field without changing the projective dimension of $I$. Moreover, we always assume that the polynomial ring $R$ is generated by a basis of the span of the entries of $M$. Since any polynomial ring containing $I$ is free over this $R$, the projective dimension of $I$ is unaffected by this choice of ambient  ring. It is easy to see that the dimension of $R$ is bounded above by the number of entries of $M$, namely $2n+2$, thus $\pd(R/I)\leq 2n+2$ follows immediately by Hilbert's Syzygy Theorem. In this section we work towards  lowering this upper bound to a {\em sharp} estimate given by $\pd(R/I)\leq 2n-2$, which we prove in Theorem \ref{linear}. 

We begin with the case when the ideal is represented by a  1-generic matrix (equivalently a matrix with no generalized zeros).

\begin{prop}\label{nozero}
Let $I$ be an ideal minimally generated by $n$ quadrics such that $I$ is contained in a prime ideal $(x,y)$ with $x,y$ independent  linear forms. Let $M$ be any  matrix representing $I$ by minors and assume that $M$ is 1-generic. Then 
\begin{enumerate}
\item $\pd (R/I)\leq n$ and
\item $I_2(M)= I:(x,y)=I:(x)=I:(y)$. 
\end{enumerate}
\end{prop}

\begin{proof}

 As shown in Lemma \ref{Cramer}, one has the containment $I_2(M)\subset I:(x,y)$. To prove that the other containment holds, we consider the possibilities for $\h(I:(x,y))$. First note that by Theorem~\ref{Eis}, $\h (I_2(M))=n$ and since $I_2(M)\subseteq I:(x,y)$, the inequality $\h (I:(x,y)) \geq n$ follows. Our first aim is to prove that in fact $\h (I:(x,y)) = n$. We  assume by way of contradiction that $\h(I:(x,y))\geq n+1$. Then, trivially, $\h((I:(x,y))+(x,y))\geq n+1$ and the hypotheses of Corollary~\ref{residual} are verified for $A=I$ and $C=(x,y)$. This yields $\h(I:(x,y))= n$, a contradiction.

We may now conclude that $\h (I:(x,y)) = n$ and we turn to computing $I:(x,y)$. Since $I_2(M)$ is by Theorem~\ref{Eis} a prime ideal  minimally generated by quadrics, neither of the linear forms $x$ and $y$ is contained in $I_2(M)$. Thus $\h (I_2(M)+(x,y)) >\h (I_2(M))$, or equivalently $\h (I_2(M)+(x,y)) \geq n+1$. Now  the hypotheses of Corollary~\ref{residual} are verified for $A=I$ and $C=(x,y)$ and an application of this Corollary yields the inequality $\pd (R/I)\leq n$, completing the proof of part (1) of our proposition. 

To prove the formula for the colon ideals in part (2), note that the ideals on both sides of the containment $I_2(M)\subseteq I:(x,y)$  have the same height $n$ and the ideal on the left is prime. This yields that in fact equality $I_2(M)= I:(x,y)$ must hold. A similar argument applies to show that the equalities $I:(x)=I:(y)=I_2(M)$ hold. 
\end{proof}

%\subsection{The case of exactly one generalized zero}

We continue our analysis of ideals $I$ generated by quadrics which are contained in a linear prime ideal of height two by considering ideals represented by matrices of the second type listed in Proposition \ref{stdzero}.

\begin{prop}\label{1zero}
Assume that $I$ is an ideal minimally generated by $n$ quadrics which is represented via minors by a matrix of the form $M=\left( \begin{matrix} x & 0 & a_{12} & a_{13} &\ldots & a_{1n} \\ y & a_{21} & a_{22} & a_{23} &\ldots & a_{2n}\\\end{matrix}\right )$, where $D=\begin{pmatrix}x &   a_{12} & \ldots & a_{1n}\\ y & a_{22} & \ldots & a_{2n}\\\end{pmatrix}$ is 1-generic. Then $\pd (R/I) \leq n$.
\end{prop}

\begin{proof}
Note that we have $I=(a_{21}x) +I'$, where $D$ represents $I'$ by minors. Furthermore, the following holds by an application of  Proposition~\ref{nozero}:
$$I:(x)=(a_{21})+I':(x)= (a_{21})+I_2(D).$$
 By Theorem~\ref{Eis},  $I_2(D)$ is a prime ideal such that $R/I_2(D)$ is Cohen-Macaulay of  codimension $n-1$. Since
the linear form $a_{21}$ is not an element of the prime ideal $I_2(D)$, it is  regular on the module $R/I_2(D)$. Together with $\pd (R/I_2(D))=n-1$, this implies that $\pd (R/(I_2(D)+(a_{21})))=n$. 
Since $I+(x)=(y)(a_{12}, a_{13},\ldots,a_{1n})+(x)$, it is easy to see that $\pd (R/(I+(x)))\leq n$. 
Since $\pd (R/(I:(x)))=\pd (R/(I_2(D)+(a_{21})))=n$ and $\pd (R/(I+(x)))\leq n$, we may now conclude by Remark \ref{pdlemma} that $\pd (R/I) \leq n$. 
\end{proof}

%\subsection{The case of at least two generalized zeros}

Next we analyze the third situation from Proposition \ref{stdzero}.

\begin{prop}\label{l1}
Assume that $I$ is an ideal  generated by $n\geq 3$ quadrics  represented by minors by a matrix $M$ which has at least two generalized zeros in the same row, i.e
is of the form
$$M=\left( \begin{matrix} x & 0 & 0 & a_{13} & \ldots & a_{1n} \\ y  &a_{21} & a_{22} & a_{23} & \ldots & a_{2n}\\\end{matrix}\right ),  $$ 
where the $a_{ij}$ are either zero or linear forms.
 Then $\pd (R/I)\leq 2n-2$.
\end{prop}

\begin{proof}
We denote by $D$ the matrix 
$$D=\left( \begin{matrix} x &  a_{13} & \ldots & a_{1n} \\ y &  a_{23} & \ldots & a_{2n}\\\end{matrix}\right ). $$
It is easy to see that $I=(xa_{21},xa_{22})+I',$
where where $D$ represents $I'$  by minors.
We conclude by Remark \ref{colon} and Proposition \ref{nozero} that
$$I: (x)= (a_{21},a_{22})+I_2(D).$$

We now set $\bar{R}=R/(a_{21},a_{22})$ and we let $\bar{D}$ denote the $2\times(n-1)$ matrix  with entries in $\bar{R}$ obtained by reducing the entries of $D$ modulo the linear forms $a_{21}$ and $a_{22}$. If $\dim \bar{R}\leq 2n-4$, we must have $\dim R\leq 2n-2$ and the desired bound  $\pd (R/I) \leq 2n-2$ follows by the Hilbert Syzygy theorem.
%Thus we are left with the case when at least one of $a_{21},a_{22}$ is not in the linear span of the entries of $D$. Without loss of generality we shall assume  $a_{21}$
%is not in the span of the entries of $D$. 

We henceforth assume $\dim \bar{R}\geq 2n-3$ (i.e. either $\dim \bar{R}= 2n-3$ or $\dim \bar{R}= 2n-2$).
If  $\bar{D}$ is 1-generic, then $\pd_{\bar{R}}(\bar{R}/I_2(\bar{D}))=n-2$ by Theorem \ref{Eis}. 
Otherwise, $\bar{D}$ has at least one generalized zero, thus,  after a linear change of coordinates in $\bar{R}$, it is equivalent to a matrix $\bar{D'}$ with $2n-3$ generic entries (each is a distinct variable) and one zero.
By a result of Boocher \cite[Theorem 4.1]{Boocher},  the minimal free resolution of $\bar{R}/I_2(\bar{D'})$ is obtained by appropriately pruning the Eagon-Northcott complex, hence  $\pd_{\bar{R}}(\bar{R}/I_2(\bar{D}))=\pd_{\bar{R}}(\bar{R}/I_2(\bar{D'}))\leq n-2$. 

To conclude, note that since $I_2(\bar{D})$ and $(a_{21},a_{22})$ use disjoint sets of variables, we have
$$
\pd (R/(I: (x))) =  \pd_R (R/((a_{21}, a_{22})+I_2(D)))  = \pd_{\bar{R}}(\bar{R}/I_2(\bar{D}))+ 2 \leq n.
 $$
It is easily seen that $\pd(R/(I+(x)))=\pd (R/(y(a_{23},\ldots,a_{2n})+(x))) \leq n-1$ and Remark \ref{pdlemma} now yields $\pd (R/I) \leq n$.
\end{proof}

\noindent Before we can analyze the next case described by Proposition \ref{stdzero},  an additional lemma is needed.

\begin{lem}\label{decomp}
 Let $I$ be an ideal generated by quadratic polynomials, contained in a linear prime ideal $(x,y)$ and represented by minors
 by a 1-generic matrix $D$. Let $a$ be a linear form  that is not in the linear span of the entries of $D$. Then
 $$(ay)+I=((y)+I)\cap((a)+I_2(D)).$$
\end{lem}

\begin{proof}
The containment $(ay)+I \subseteq ((y)+I)\cap((a)+I_2(D))$ is obvious. As for the opposite containment, let $f\in (y+I)\cap((a)+I_2(D))$. Since 
$f\in ((y)+I)$, we have $f=yh+g$, with $g\in I$ and $h\in R$. As $g\in I \subset (a)+I_2(D)$, the assertion $f\in (a)+I_2(D)$ is equivalent to $yh\in (a)+I_2(D)$.

By Theorem~\ref{Eis}, $I_2(D)$ is a prime ideal. Since  $a$ is a linear form which is is not in the span of the entries of $D$,   the ideal $(a)+I_2(D)$ is in turn prime. Now $yh\in (a)+I_2(D)$ implies $y\in  (a)+I_2(D)$ or $h\in  (a)+I_2(D)$. Since the only linear forms in $ (a)+I_2(D)$ are scalar multiples of $a$, which is
 not a multiple of $y$ by hypothesis, the first alternative is impossible. It remains that $h\in  (a)+I_2(D)$, thus finishing the proof.
\end{proof}

We continue our analysis with the fourth among the cases listed in Proposition \ref{stdzero}.

\begin{prop}\label{l2}
Let $I$ be an ideal  generated by $n\geq 3$ quadrics  and represented by minors by a matrix $M$  of the form
$$M=\left( \begin{matrix} x & 0 & a_{12} & a_{13} & \ldots & a_{1n} \\ y & a_{21} & 0 & a_{23} & \ldots & a_{2n}\\\end{matrix}\right ), $$ where the submatrix $D=\left( \begin{matrix} x &  a_{13} & \ldots & a_{1n} \\ y &  a_{23} & \ldots & a_{2n}\\\end{matrix}\right ) $ is 1-generic. Then $\pd (R/I)\leq 2n-2$.
\end{prop}

\begin{proof}
Clearly the ideal $I$ has the form $I=(a_{12}y,a_{21}x)+I'$, where $I'$ is the ideal represented by minors by the 1-generic matrix $D$. Should both $a_{21}$ and $a_{12}$ be in the span of the entries of $D$, then the conclusion would easily follow from Hilbert's Syzygy Theorem. Therefore, without loss of generality, we may assume that $a_{12}$ is not in the span of the entries of $D$.   Furthermore, if $a_{21}$ is in the span of $\{y,a_{23},\ldots,a_{2n}\}$, then the second column of $M$ can be replaced by ideal-preserving matrix operations by a column having a zero entry in the bottom, and then the ideal $I$ can be represented by minors by a matrix that has already been analyzed in Proposition \ref{l1}. Henceforth, we assume that $a_{21}$ is not in the span of $\{y,a_{23},\ldots,a_{2n}\}$. 

It is easy to see that  $\pd(R/(I+(x)))=\pd(R/(x,ya_{12},ya_{13},\ldots,ya_{1n}))\leq n$. We now proceed to compute $I:(x)$ using Remark \ref{colon} and Lemma \ref{decomp}: 
\[
\begin{split}
I:(x) &= (a_{21})+(a_{12}y+I'):(x)=(a_{21})+[((y)+I')\cap((a_{12})+I_2(D))]:(x)\\
 &=  (a_{21}) +[(y,a_{23}x,\ldots,a_{2n}x):(x)]\cap[((a_{12})+I_2(D)):(x))]\\
 &=  (a_{21}) +(y,a_{23},\ldots,a_{2n})\cap((a_{12})+I_2(D)).\\
\end{split}
\]
The last equality follows because the ideal $(a_{12})+I_2(D)$ is a prime ideal not containing the any linear form in the span of $x$. 
We now estimate the projective dimension of each term appearing above. First, by Theorem \ref{Eis} we have $\pd(R/I_2(D))=n-2$  and since $a_{12}$ is a regular form on $I_2(D)$, we obtain that $\pd(R/(a_{12}+I_2(D))=n-1$. Furthermore, we have  that $\pd(R/(y,a_{23},\ldots,a_{2n}))=n-1$ and 
$$\pd(R/((y,a_{23},\ldots,a_{2n})+(a_{12}+I_2(D))))=\pd(R/(y,a_{23},\ldots,a_{2n},a_{12}))\leq n.$$ Using a standard short exact sequence and the three previous inequalities we deduce: $$\pd(R/((y,a_{23},\ldots,a_{2n})\cap((a_{12})+I_2(D))))\leq n-1.$$

\noindent Since $(y, a_{23}, ... , a_{2n}) \cap ((a_{12}) + I_2(D)))$ is
the intersection of two prime ideals and since $a_{21}\notin (y, a_{23}, ... , a_{2n})$, it is easily seen that $a_{21}$ is a non-zero divisor on $R/((y, a_{23}, ... , a_{2n}) \cap ((a_{12}) + I_2(D)))$, if and only if $a_{21}$ is not in the span of $a_{12}$. Thus, if  $a_{21}$ is not in the span of $a_{12}$, then we deduce from the previous inequality that
$$\pd(R/((a_{21})+(y,a_{23},\ldots,a_{2n})\cap((a_{12})+I_2(D)))\leq n.$$ Since $\pd(R/(I+(x))=\pd(R/(x,ya_{12},ya_{13},\ldots,ya_{1n}))\leq n$, it follows that $\pd(R/I)\leq n$.

 % THIS WAS SHORTENED We claim that $a_{21}$ a non-zero divisor on  $R/((y,a_{23},\ldots,a_{2n})\cap((a_{12})+I_2(D)))$, unless $(a_{21})=(a_{12})$. Indeed, suppose that $(a_{21})\neq (a_{12})$ and $a_{21}F\in (y,a_{23},\ldots,a_{2n})\cap((a_{12})+I_2(D))$ for some polynomial $F$. Note that this is the intersection of two prime ideals and since $a_{21}$ is not an element of either ideal, it follows that $F\in(y,a_{23},\ldots,a_{2n})\cap((a_{12})+I_2(D))$.  We may now deduce that $$\pd(R/((a_{21})+(y,a_{23},\ldots,a_{2n})\cap((a_{12})+I_2(D)))\leq n.$$ Since $\pd(R/(I+(x))=\pd(R/(x,ya_{12},ya_{13},\ldots,ya_{1n}))\leq n$, it follows that $\pd(R/I)\leq n$.

 The remaining case is when $a_{12}$ is a scalar multiple of $a_{21}$. Recall that $a_{12}$ is not in the span of the entries of $D$, which allows us to compute 
$I:(a_{12})=(x,y)+I':(a_{12})=(x,y)+I'=(x,y) \mbox{ and } I+(a_{12})=I'+(a_{12}).$
By Proposition \ref{nozero}, we have that $\pd(R/I')\leq n+1$, since $D$ is 1-generic. Since $a_{12}$ is a non-zerodivisor on $R/I'$, it follows that $\pd(R/((a_{12}+I'))\leq n+2$; therefore $\pd(R/I)\leq n+1$ in this case.
\end{proof}

%% ANOTHER TRY %%%%
%We now show that $$(y,a_{23},\ldots,a_{2n})\cap((a_{12})+I_2(D))=I_2(C), \mbox{ where } 
%C=\left( \begin{matrix} x &  a_{12} &a_{13} & \ldots & a_{1n} \\ y &  0 & a_{23} & \ldots & a_{2n}\\\end{matrix}\right ). $$
%Indeed, since $I_2(D)\subseteq (y,a_{23},\ldots,a_{2n})$, we have that $(y,a_{23},\ldots,a_{2n})\cap((a_{12})+I_2(D))=(y,a_{23},\ldots,a_{2n})\cap(a_{12})+I_2(D)= (a_{12}y,a_{12}a_{23},\ldots,a_{12}a_{2n})+I_2(D)=I_2(C)$. It follows from our computation above that $I:(x)=(a_{21})+I_2(C)$.
%To estimate the projective dimension of this colon ideal we begin by defining $\bar{R}=R/(a_{21})$ and letting $\bar{C}$ and $\bar{D}$ denote the matrix obtained by reducing the entries of $C$ or $D$ respectively modulo $(a_{21})$. Then
%$$\pd_R(R/((a_{21})+I_2(C))=\pd_{\bar{R}}(I_2(\bar{C}))+1.$$
%Note that $\bar{C}$ has at most $2n-1$ non-zero entries. If the non-zero entries of $\bar{C}$ form a linearly independent set, then a result of Boocher  in \cite[]{Boocher} implies that $\pd_{\bar{R}}(I_2(\bar{C}))\leq n-1$. Otherwise we have that the non-zero entries of $\bar{C}$ form a linearly dependent set. We may assume that the dimension of this set is $2n-2$, otherwise an applicatieon of Hilbert's Syzygy Theorem finishes the proof. If $\bar{D}$ is not 1-generic, we may find a matrix equivalent to $C$ having exactly 2 zero entries and $2n-2$ non-zero entries which form a linearly independent set.  In this case Boocher's result yields that $\pd_{\bar{R}}(I_2(\bar{C}))\leq n-1$ again. It remains to consider the case when  $\bar{D}$ is 1-generic.

Finally, we consider the last  case described in Proposition \ref{stdzero}. To handle this case, we need to assume an inductive bound on the projective dimension of ideals of quadrics with fewer generators. We shall use the Proposition below in an inductive argument to compute the overall bound given in Proposition \ref{linear}.

\begin{prop}\label{l3}
Let $I$ be an ideal  generated by $n\geq 3$ quadrics  and represented by minors by a matrix $M$  of the form
$$M=\left( \begin{matrix} x & 0 & a_{12} & a_{13} & a_{14} & \ldots & a_{1n} \\ y & a_{21} & 0 & \lambda a_{13} & a_{24} & \ldots & a_{2n}\\\end{matrix}\right ), $$ where the $a_{ij}$ are either zero or linear forms.
Assume additionally that any ideal $J$ of height two generated by $n-1$ quadrics in a polynomial ring $S$ is known to have $\pd (S/J)\leq 2n-4$.
Then $\pd (R/I)\leq 2n-2$. 
\end{prop}

\begin{proof}
Note that, if the non-zero entries of $M$ form a linearly dependent set, then the conclusion follows immediately by the Hilbert Syzygy Theorem. We henceforth assume that the $2n-1$ non-zero entries of $M$ form a linearly independent set.
It can be read off $M$ that the ideal $I$ has the form
$$I=(a_{12}y,a_{21}x, a_{13}(\lambda x- y))+I',$$
where $I'$ is an ideal generated by quadrics represented by minors by the  matrix given by the first column and last $n-2$ columns of $M$. Next we compute 
$$I:(a_{12})=(y)+((a_{21}x, a_{13}(\lambda x- y))+I'):(a_{12})=(y)+(a_{21}x, a_{13}(\lambda x- y))+I',$$
where the first equality follows from Remark \ref{colon} and the second from  knowing that the ideal $(a_{21}x, a_{13}(\lambda x- y)+I'$ is written in terms of a set of variables disjoint from $a_{12}$. Since  the previous computation gives
$I:(a_{12})=(y)+(a_{21},a_{13}, a_{24},\ldots, a_{2n})x$, 
we may conclude that $\pd (R/(I:a_{12}))=n$.

Next we compute $I+(a_{12})$. Since $a_{12}$ only appears in one minimal generator of $I$, one has
$I+(a_{12})=J+(a_{12})$, where $J$ in the ideal minimally generated by the $n-1$ quadratic generators of $I$ not involving $a_{12}$. Consider the polynomial
 ring $S=K[x, y, a_{21}, a_{13}, a_{ij}]$, (for $1\leq i \leq 2$ and $4\leq j\leq n$)   and view $J$ as an ideal of  $S$. The hypothesis on ideals generated by $n-1$ quadrics yields
 $\pd (S/J)\leq 2n-4$. Since $a_{12}$ is a variable not appearing in $S$,
 we have $\pd_R (R/((a_{12})+J))=\pd_S (S/J)+1 \leq 2n-3$.
 Finally, we use Remark \ref{pdlemma} to estimate $\pd (R/I) \leq \max \{n,2n-3\} \leq 2n-2$ for all $n\geq 3$.   
\end{proof}

%\subsection{The overall bound}

To complete the analysis, we  combine Propositions \ref{nozero}, \ref{1zero}, \ref{l1} and \ref{l2}  to obtain:

\begin{thm}\label{linear}
Assume that $I$ is a height two ideal minimally generated by $n$ quadrics which has a  linear minimal associated prime. Then
$\pd (R/I) \leq 2n-2$.
\end{thm}

\begin{proof}
We proceed by induction on $n$, noting that the assumption on the height of $I$ forces $n\geq 2$. The base case, $n=2$ is easily seen to be  a complete intersection thus having projective dimension two.

Let $n\geq 3$. Using the inductive hypothesis, we may  assume that any ideal $J$ generated by $n-1$ quadrics, which has a minimal associated prime ideal  generated by two linear forms has $\pd (S/J)\leq 2n-4$,  where $S$ is an ambient polynomial ring for $J$. Now we consider an ideal $I$ generated by $n$ quadrics having a linear minimal  prime ideal of height two. 

We finish the proof by analyzing the various possibilities for matrices representing $I$ by minors listed in  Proposition \ref{stdzero}:
\begin{enumerate}
\item if $M$ has no generalized zeros, we apply Propositions~\ref{nozero} to obtain $\pd (R/I)\leq n+1$;
\item if $M$  is of one of the types stated in  Propositions~\ref{1zero}, \ref{l1} or \ref{l2}, then we obtain  $\pd (R/I)\leq 2n-2$;
\item if $M$ is of the type stated in Proposition \ref{l3} we rely on the inductive hypothesis to obtain again $\pd (R/I)\leq 2n-2$.
\end{enumerate}
\end{proof}

\subsection{Better bounds}

In Theorem~\ref{linear} we have successfully bounded the projective dimension of ideals generated by $n\geq 2$ quadrics under the assumption that the ideal is contained in a linear prime ideal of height two. We now improve on this bound under the additional condition that the multiplicity of the ideal at this linear prime is at least two.  We recall that, by Lemma \ref{Pmult}, the local multiplicity of an ideal of height two generated by quadrics at such a  prime is at most three.

\begin{prop}\label{1;n}
Let $R$ be a polynomial ring; let $I$ be an ideal generated by $n\geq 2$ quadrics, contained in a height two linear prime ideal $(x,y)$ and assume $e(R_{(x,y)}/I_{(x,y)})\geq 2$. Then $\pd (R/I) \leq n+2$.
\end{prop}
\begin{proof} 
Let $M$ be any matrix representing $I$ by minors and let $A$ be the  submatrix represen-ting $I$ by coefficients. Since $e(R_{(x,y)}/I_{(x,y)})\geq 2$,  we must have that $I_2(M)R_{(x,y)}\neq R_{(x,y)}$; otherwise, the containment $(x,y)R_{(x,y)}\subseteq IR_{(x,y)}$ would hold by Lemma \ref{Cramer}, leading to $e(R_{(x,y)}/I_{(x,y)})=1$. The fact that $I_2(M)R_{(x,y)}\neq R_{(x,y)}$ is equivalent to $I_2(M)\subseteq (x,y)$, which in particular implies that the containment $I_2(A)\subseteq(x,y)$ holds. Now consider  the matrix $\bar{A}$ obtained by taking the image of $A$ under the canonical projection to the polynomial factor ring  $\bar{R}=R/(x,y)$. Since $I_2(A)\subseteq(x,y)$, it follows that $I_2(\bar{A})=(0)$ in $\bar{R}$. An application of Theorem \ref{EHU} for $\bar{A}$, with $r=1$ and $i=1$, shows that the entries of $\bar{A}$ can be written in terms of at most $n$ independent linear forms. The entries of the matrix $A$ and hence also the ideal $I$ can then be written in terms of at most $n+2$ independent linear forms. An application of the Hilbert Syzygy Theorem finishes the proof.
\end{proof}

 Proposition \ref{xq} on ideals contained in a prime of multiplicity two can be imitated in the context of bounding the projective dimension of ideals contained in linear primes as well, with the additional hypothesis that there be at least two distinct such linear minimal primes.

\begin{prop}\label{2lin}
Let $I$ be an ideal generated by $n$ quadratic polynomials which is contained in two distinct linear prime ideals of height 2. Then $\pd (R/I) \leq n$.
\end{prop}

\begin{proof}
Let the two linear primes be $(x,y)$ and $(z,w)$. If $\h (x,y,z,w)=4$, then the intersection of these two ideals is generated by quadrics and since the minimal generators of $I$ must be linear combinations of these quadrics,the variables $x,y,z,w$ are the only ones required to express  $I$, hence $\pd (R/I) \leq 4$. Otherwise if 
$\h (x,y,z,w)<4$, we may assume without loss of generality that $x=w$. Now $I\subset (x,y)\cap(x,z)=(x,yz)$. Setting $q=yz$, one may apply the technique in the proof of Proposition~\ref{xq} to see that $\pd (R/I)\leq n$.
\end{proof}

\section{Final questions}

Throughout this paper two key invariants were used to obtain our upper bound on $\pd (R/I)$: the number of minimal generators of $I$ and the additional assumption that $\h I=2$. In this section, we formulate two questions aimed towards establishing whether these invariants suffice in general for bounding the projective dimension of ideals generated by quadrics. The main result of this paper answers both question below in the affirmative for the case $h=2$. Together with a result of Eisenbud-Huneke (cf. \cite[Theorem 3.1]{MS}) for the case $n=3$ and easy arguments in the remaining cases, this verifies that the two questions have an affirmative answer whenever $n\leq3$ or $h\leq 2$. 

\begin{qst}\label{qst1}
Let $R$ be a polynomial ring and let $I$ be an ideal of $R$ generated by $n$ quadrics and having  $\h I=h$. Is there an upper bound for $\pd(R/I)$ expressed only in terms of $n$ and $h$?
\end{qst}

Based on computational evidence and on the bound given by our main result in the case of height two ideals of quadrics, we refine the previous question by suggesting a specific candidate for an upper bound on the projective dimension of ideals of quadrics.

\begin{qst}\label{qst2}
Let $R$ be a polynomial ring and let $I$ be an ideal of $R$ generated by $n$ quadrics and having  $\h I=h$. Is it true that $\pd (R/I)\leq h(n-h+1)$?
\end{qst}

We note that there is in the a family of ideals of quadrics which achieves the equality in the bound proposed in Question \ref{qst2}, namely the ideals  $I_{m,n,2}$ described in  \cite{McCullough}. Therefore, if the answer to Question \ref{qst1} is affirmative, any such bound must be at least as large as the bound we suggest in Question \ref{qst2}.
However, as soon as the generators are allowed to have arbitrary degrees, one cannot expect any upper bound on the
projective dimension to have asymptotic growth bounded the expression in  Question \ref{qst2} or even by the product of height and minimal number of generators. To see this, consider again, this time for $d\geq 3$, the family of ideals $I_{m,n,d}$ described in  \cite{McCullough}. These ideals have $\h I_{m,n,d}\leq m$ and are minimally generated by $m+n$ homogeneous polynomials of degree $d$, while their projective dimension is given by $\pd (R/I_{m,n,d})=m+n\frac{(m+d-2)!}{(m-1)!(d-1)!}$.

\section*{Acknowledgements}
This work was completed while all the authors were members in the MSRI program on Commutative Algebra. We thank  the organizers of the program  and  MSRI for making this possible. The first author was supported by NSF grant number 1259142.

%\newpage

\bibliographystyle{amsplain}
\bibliography{bibht2}

\end{document}